\documentclass[12pt]{article}
\usepackage{latexsym}
\newtheorem{defi}{Definition}
\newtheorem{theorem}{Theorem}
\title{Extreme Palindromes}
\author{Kathy Q. Ji\\
Center for Combinatorics, LPMC, Nankai University\\
Tianjin 300071, P. R. China\\
\texttt{<ji@nankai.edu.cn>}
\and
Herbert S. Wilf\\
Department of Mathematics, University of Pennsylvania\\
Philadelphia, PA 19104-6395\\
\texttt{<wilf@math.upenn.edu>}
}
\begin{document}
\maketitle
\begin{abstract}
A recursively palindromic (RP) word is one that is a palindrome and whose left half-word and right half-word are
 each RP. Thus ABACABA is, and MADAM is not, an RP word. We count RP words of given length over a finite alphabet
 and RP compositions of an integer. We use the same method to determine the parity of the Catalan numbers.
\end{abstract}
\section{Introduction}
Suppose we have a collection of tuples and we are interested in the
parity of the size of the collection. One obvious method of reducing
the problem would be to define an involution on the tuples, and then
study the parity of the set of its fixed points. One way to define
such an involution would be to associate each tuple with its
reverse, and discard the pair. That would leave us with only the
 palindromic tuples to study. But we can go further.

Suppose $(a_1,a_2,\dots,a_{2k})$ is one of the remaining tuples of
even length. Even though it is palindromic, it might be that its
left half $(a_1,\dots,a_k)$ is not, in which case its right half is
also not palindromic. In that case we could pair this tuple with the
one obtained by reversing its left half and reversing its right
half, and similarly, recursively, at all levels. If
$(a_1,a_2,\dots,a_{2k+1})$ is one of the remaining tuples of odd
length, then even though it is palindromic, it might be that its
left half $(a_1,\dots,a_k)$ is not, etc. as before.

The end result of this pairing and discarding process, i.e., the set
of unpaired tuples, would be the collection of tuples that are
\textit{recursively palindromic} (RP), assuming that the collection
of tuples is closed under the various reversal operations.

Thus instead of studying the parity of the number of all tuples in
the given collection, it would suffice to study the parity of the
collection of RP tuples.

\begin{defi}
A tuple is recursively palindromic if it is empty, or it is a palindrome and its
left half and its right half are recursively palindromic.
\end{defi}

The word ABACABA is RP but the word MADAM is not.

\section{The pairing}
In general, we are given a collection ${\cal C}$ of words over some
alphabet. We suppose that for every $w\in {\cal C}$, if $\sigma(w)$
is some rearrangement of the letters of $w$, then $\sigma(w)\in{\cal
C}$ also. Then we associate with each word $w\in {\cal C}$ a binary
tree $T(w)$, as follows.

By the \textit{left half} of a word $w=(a_1a_2\dots a_{2k})$ we mean
the word $w_L=(a_1a_2\dots a_{k})$, while the left half of
$w=(a_1a_2\dots a_{2k+1})$ will also be $w_L=(a_1a_2\dots a_{k})$,
and similarly for the right half $w_R$ of a word. Then a labeled
binary tree $T(w)$, associated with the word $w$, is constructed
recurrently as follows: for the word $w=(a_1a_2\dots a_{2k})$, label
the root of $T$ by $\emptyset$, while for the word $w=(a_1a_2\dots
a_{2k+1})$, label the root by $a_{k+1}$. Further, the left subtree
at the root of $T(w)$ is $T(w_L)$ and the right subtree at the root
of $T(w)$ is $T(w_R)$.

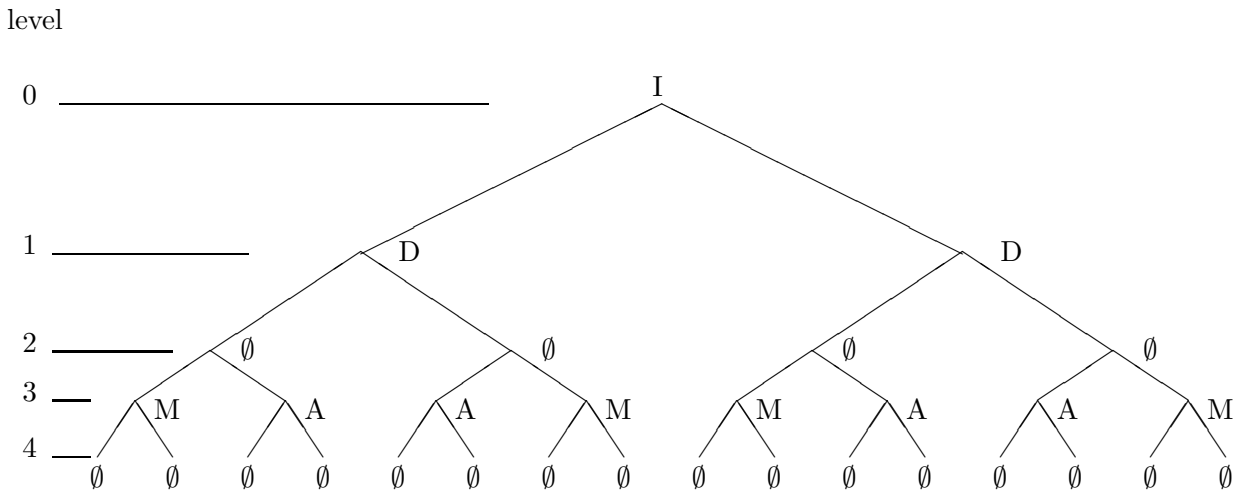
\begin{figure}[h]
\begin{center}
\begin{picture}(330,200)
\setlength{\unitlength}{1mm}
%------------------------------------------------------------------
\put(-50,10){\put(40,0){\line(2,3){5}}\put(50,0){\line(-2,3){5}}\put(70,0)
{\line(-2,3){5}}\put(60,0){\line(2,3){5}}
\put(80,0){\line(2,3){5}}\put(59,-4){\small{$\emptyset$}}\put(49,-4)
{\small{$\emptyset$}}
\put(69,-4){\small{$\emptyset$}}\put(79,-4){\small{$\emptyset$}}
\put(90,0){\line(-2,3){5}}\put(39,-4){\small{$\emptyset$}}
\put(89,-4){\small{$\emptyset$}}\put(99,-4){\small{$\emptyset$}}
\put(109,-4){\small{$\emptyset$}}
\put(119,-4){\small{$\emptyset$}}\put(129,-4){\small{$\emptyset$}}
\put(139,-4){\small{$\emptyset$}}\put(149,-4){\small{$\emptyset$}}
\put(159,-4){\small{$\emptyset$}}\put(169,-4){\small{$\emptyset$}}
\put(179,-4){\small{$\emptyset$}}\put(189,-4){\small{$\emptyset$}}
%\put(89,-4){\small{A}}
\put(100,0){\line(2,3){5}}
\put(110,0){\line(-2,3){5}}\put(120,0){\line(2,3){5}}
\put(130,0){\line(-2,3){5}}\put(140,0){\line(2,3){5}}
\put(150,0){\line(-2,3){5}}\put(160,0){\line(2,3){5}}
\put(170,0){\line(-2,3){5}}\put(180,0){\line(2,3){5}}
\put(190,0){\line(-2,3){5}} \put(127.5,5){\small{M}}
\put(167.5,5){\small{A}}
\put(187.5,5){\small{M}}\put(107.5,5){\small{M}}
\put(87.5,5){\small{A}}\put(67.5,5){\small{A}}
\put(47.5,5){\small{M}}
\put(65,7.5){\line(-3,2){10}}\put(45,7.5){\line(3,2){10}}\put(85,7.5){\line(3,2){10}}
\put(105,7.5){\line(-3,2){10}}\put(125,7.5){\line(3,2){10}}
\put(145,7.5){\line(-3,2){10}}\put(165,7.5){\line(3,2){10}}
\put(185,7.5){\line(-3,2){10}}\put(99,13){\small{$\emptyset$}}
\put(139,13){\small{$\emptyset$}}\put(179,13){\small{$\emptyset$}}\put(59,13){\small{$\emptyset$}}
\put(147.5,5){\small{A}}\put(55,14){\line(3,2){20}}
\put(95,14){\line(-3,2){20}}\put(135,14){\line(3,2){20}}
\put(175,14){\line(-3,2){20}}\put(160,26){\small{D}}
\put(80,26){\small{D}}
\put(75,27){\line(2,1){40}}\put(155,27){\line(-2,1){40}}
\put(112.5,48){\small{ I}}
\put(39,0){\line(-1,0){5}}\put(30,0){\small{4}}
\put(39,7.5){\line(-1,0){5}}\put(30,7.5){\small{3}}\put(50,14){\line(-1,0){16}}
\put(60,27){\line(-1,0){26}}\put(92,47){\line(-1,0){57}}
\put(30,14){\small{2}}\put(30,27){\small{1}}\put(30,47){\small{0}}\put(28,57)
{\small{level}}}
\end{picture}
\end{center}
\caption{The binary tree of height 5 corresponding to
``MADAMIMADAM''} \label{fg:bt1}
\end{figure}

  The binary tree $T(\mathrm{MADAMIMADAM})$, for example, is shown in
 Fig. \ref{fg:bt1}. If the tree $T(w)$ is given then we can recover
 the word $w$ by concatenating the labels at the nodes of $T(w)$
  when they are traversed in inorder (left-root-right).

Here we remark that the binary tree of an RP word is a binary tree
 such that all nodes on the same level have the same label.

We now use the tree $T(w)$ of a word to define an involution on the
collection ${\cal C}$. Given a word $w$ which is not RP, to find the
word $w'$ that is paired with $w$, proceed as follows. Begin with
the root of $T(w)$, and go down the levels of $T$ until for the
first time reaching a level such that the labels of the nodes on
level $L$ are not all the same. Then interchange the right and left
subtrees at every node on level $L-1$, to obtain $T(w')$. Then we
obtain $w'$ by inorder traversal of $T(w').$ It's obvious that this
map is an involution.

In the case of MADAMIMADAM, for example, we would exchange the left
and right subtrees at each node on level 2 of the tree in Fig.
\ref{fg:bt1} because the labels on level 3 are not all the same. By
visiting the nodes of the new tree in inorder, we would find that
MADAMIMADAM has been paired with AMDMAIAMDMA by this mapping.

\section{RP words}
Suppose we have an alphabet of $K$ letters. Of the $K^n$ possible words of length $n$,
 how many are RP words?
\begin{theorem}
There are exactly $K^{\alpha(n)}$ RP words of length $n$ over an alphabet of $K$ letters,
 where $\alpha(n)$ denotes the sum of the binary digits of $n$.
\end{theorem}

The easy proof of this theorem is by recurrence. Let $f(n)$ be the required number of
words. Then evidently $f(2n)=f(n)$, for $n\ge 1$, and $f(2n+1)=Kf(n)$, for $n\ge 0$,
 with $f(0)=1$. The function $K^{\alpha(n)}$ satisfies the same recurrences with the
 same initial value. $\Box$

For example, among the 128 binary words of length 7, there are 8 RP words, viz.
\[ 0000000,\,0001000,\,0100010,\,0101010,\,1010101,\,1011101,\,1110111,\,1111111.\]
\subsection{Bijective proof}
Given the binary representation of $n$. Consider the collection of
all sequences $\pi$ which can be obtained from that binary
representation by replacing each 1  by some letter in the given
alphabet of $K$ letters. Evidently $K^{\alpha(n)}$ is the number of
 such sequences. If $\pi$ is such a sequence,
 let its length be $L.$ We will first build a
bijection between these sequences $\pi$ and complete binary trees
$CT(\pi)$ of height $L$ whose vertices are labeled with letters from
our alphabet.

If $\pi=\emptyset$, then $CT(\pi)=\emptyset.$ If $\pi\neq
\emptyset$, then for each $i$, position $i$ in the sequence $\pi$
will correspond to level $i$ of a complete binary tree. Furthermore
if the letter in position $i$ of the sequence is $t$ , then all
nodes on level $i$ of that complete binary tree are labeled by $t$.
Hence, for example, the sequence $\pi=(A0B0C)$ will correspond to
the following complete binary tree of height 5:

\begin{figure}[h]
\begin{center}
\begin{picture}(330,200)
\setlength{\unitlength}{1mm}
%------------------------------------------------------------------
\put(-50,10){\put(40,0){\line(2,3){5}}\put(50,0){\line(-2,3){5}}\put(70,0)
{\line(-2,3){5}}\put(60,0){\line(2,3){5}}
\put(80,0){\line(2,3){5}}\put(59,-4){\small{A}}\put(49,-4){\small{A}}
\put(69,-4){\small{A}}\put(79,-4){\small{A}}
\put(90,0){\line(-2,3){5}}\put(39,-4){\small{A}}
\put(89,-4){\small{A}}\put(99,-4){\small{A}}\put(109,-4){\small{A}}
\put(119,-4){\small{A}}\put(129,-4){\small{A}}
\put(139,-4){\small{A}}\put(149,-4){\small{A}}
\put(159,-4){\small{A}}\put(169,-4){\small{A}}
\put(179,-4){\small{A}}\put(189,-4){\small{A}}
\put(89,-4){\small{A}}\put(100,0){\line(2,3){5}}
\put(110,0){\line(-2,3){5}}\put(120,0){\line(2,3){5}}
\put(130,0){\line(-2,3){5}}\put(140,0){\line(2,3){5}}
\put(150,0){\line(-2,3){5}}\put(160,0){\line(2,3){5}}
\put(170,0){\line(-2,3){5}}\put(180,0){\line(2,3){5}}
\put(190,0){\line(-2,3){5}} \put(124.5,3){\small{0}}
\put(144.5,3){\small{0}}\put(164.5,3){\small{0}}
\put(184.5,3){\small{0}}\put(104.5,3){\small{0}}
\put(84.5,3){\small{0}}\put(64.5,3){\small{0}}
\put(44.5,3){\small{0}}
\put(65,7.5){\line(-3,2){10}}\put(45,7.5){\line(3,2){10}}\put(85,7.5){\line(3,2){10}}
\put(105,7.5){\line(-3,2){10}}\put(125,7.5){\line(3,2){10}}
\put(145,7.5){\line(-3,2){10}}\put(165,7.5){\line(3,2){10}}
\put(185,7.5){\line(-3,2){10}}\put(94,10){\small{B}}
\put(134,10){\small{B}}\put(174,10){\small{B}}\put(54,10){\small{B}}
\put(55,14){\line(3,2){20}}
\put(95,14){\line(-3,2){20}}\put(135,14){\line(3,2){20}}
\put(175,14){\line(-3,2){20}}\put(155,23){\small{0}}
\put(75,23){\small{0}}
\put(75,27){\line(2,1){40}}\put(155,27){\line(-2,1){40}}\put(115,42){\small{C}}
\put(39,0){\line(-1,0){5}}\put(30,0){\small{4}}
\put(39,7.5){\line(-1,0){5}}\put(30,7.5){\small{3}}\put(50,14){\line(-1,0){16}}
\put(60,27){\line(-1,0){26}}\put(100,47){\line(-1,0){66}}
\put(30,14){\small{2}}\put(30,27){\small{1}}\put(30,47){\small{0}}\put(28,57)
{\small{level}}}
\end{picture}
\end{center}
\caption{A complete binary tree of height 5 corresponding to an RP
word.}\label{twintree}
\end{figure}
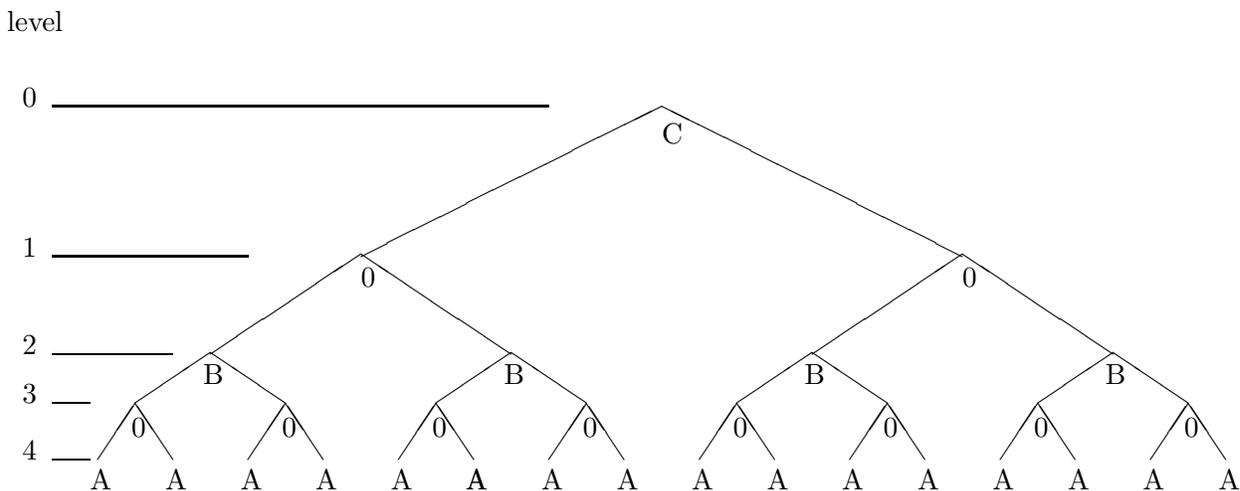
Note that the number of \textit{nodes} labeled by the letters in
such a complete binary tree is equal to $n$, while the number of
\textit{levels} labeled by a letter is equal to the sum of the
binary digits of $n$, namely $\alpha(n)$. Then, reading all letters
on such a complete binary tree by  \textit{inorder traversal}
(left-root-right) will give an RP word of $n$ letters over the
alphabet $\{A,B,C,\ldots,\,K\}$. Thus, e.g., the sequence $(A0B0C)$
corresponds to the RP word $AABAAAABAACAABAAAABAA$ of $21$ letters.

%The bijection is easily described in general. For an RP word
%$\omega$ of $n$ letters over the alphabet $\{A,B,C,\ldots,\,K\}$,
%we will construct a complete binary tree $CT(\omega)$ labeled by
%$0,A,\ldots,K.$ If $\omega=\emptyset$, then
%$CT(\omega)=\emptyset$. If $\omega\neq \emptyset$ and $n$ is odd,
%let the letter in the position $\frac{n-1}{2}$ be $t$. Thus
%$\omega=\omega't\omega'$ where $\omega'$ is an RP words of
%$\frac{n-1}{2}$ letters. Now let $t$ be the root of $CT(\omega)$
%and let $CT(\omega')$  be the left  subtree obtained by removing
%$t$ as well as  right subtree.  If $n$ is even, then
%$\omega=\omega'\omega'$ where $\omega'$ is an RP words of
%$\frac{n}{2}$ letters. Now let $0$ be the root of $CT(\omega)$,
%and let $CT(\omega')$ be the left  subtree as well as  right
%subtree. This yields an inductive definition of $CT(\omega).$
It's easy to see in general that the resulting $CT(\omega)$ is a
complete binary tree, and furthermore all nodes on the same level
have the same label. Hence it will correspond to a sequence counted
by $K^{\alpha(n)}.$ $\Box$
\section{RP compositions}
By a composition of $n$ we will mean an ordered representation of $n$
as a sum of positive integers, called the \textit{parts} of the composition.
 For instance, $6=1+2+1+2$ is a composition of 6 into four parts. There are
  $2^{n-1}$ compositions of $n$. Notice that, for example, $2+1+2+6+2+1+2$ is
   a RP composition of $n=16$.  We ask for the number of compositions of $n$
    that are RP.

Let $f(n)$ be the required number. Now an RP composition of $2n$ is either of
the form
\begin{enumerate}
\item $ww$, where $w$ is an RP composition of $n$, or
\item $w(2j)w$, where $j>0$ and $w$ is an RP composition of $n-j$.
\end{enumerate}
Thus $f(2n)=\sum_{j=0}^nf(n-j)$, and by subtraction, $f(2n)=f(2n-2)+f(n)$.

An RP composition of $2n+1$ is of the form $w(2j+1)w$, for some $j\ge 0$,
where $w$ is an RP composition of $n-j$. Thus $f(2n+1)=\sum_{j=0}^nf(n-j)=f(2n)$.
Hence the counting function $f(n)$ satisfies
\[f(2n+1)=f(2n)\quad(n\ge 0);\qquad f(2n)=f(2n-2)+f(n)\quad(n\ge 1);\qquad f(0)=1.\]
These recurrences and the initial value are identical with those satisfied by the
 sequence $b(n)$, the number of partitions of $n$ into powers of 2, and we have
\begin{theorem}
The number of RP compositions of $n$ is equal to the number of
partitions of $n$ into powers of 2 {\rm(}``binary
partitions''{\rm)}.
\end{theorem}
\section{Bijective proof}
For a partition $\lambda$ of $n$ into powers of $2$, suppose that
its largest part is $2^l$. Similarly to the bijection on RP words,
we will first build a bijection between the partitions $\lambda$ of
$n$ into powers of $2$ and complete binary trees $CT(\lambda)$ of
height $l+1$ whose vertices are labeled with nonnegative integers.

If $\lambda=\emptyset$, then $CT(\lambda)=\emptyset.$ If
$\lambda\neq \emptyset$, then for each $i=0,1,\dots,l$, the part
$2^i$
 of $\lambda$  will correspond to the level
$i$ of the complete binary tree, which has $2^i$ nodes.
Furthermore if the part $2^i$ of $\lambda$ has multiplicity $m_i$, then all
nodes on level $i$ of the complete binary tree are labeled by
$m_i$. Hence the partition $\lambda=(16, 4, 4,4,4,1,1,1,1,1)$ of
37 will correspond to the complete binary tree of
height 5 that is shown in Figure \ref{twintree}.

\begin{figure}[h]
\begin{center}
\begin{picture}(330,200)
\setlength{\unitlength}{1mm}
%------------------------------------------------------------------
\put(-50,10){\put(40,0){\line(2,3){5}}\put(50,0){\line(-2,3){5}}\put(70,0)
{\line(-2,3){5}}\put(60,0){\line(2,3){5}}
\put(80,0){\line(2,3){5}}\put(59,-4){\small{1}}\put(49,-4){\small{1}}
\put(69,-4){\small{1}}\put(79,-4){\small{1}}
\put(90,0){\line(-2,3){5}}\put(39,-4){\small{1}}
\put(89,-4){\small{1}}\put(99,-4){\small{1}}\put(109,-4){\small{1}}
\put(119,-4){\small{1}}\put(129,-4){\small{1}}
\put(139,-4){\small{1}}\put(149,-4){\small{1}}
\put(159,-4){\small{1}}\put(169,-4){\small{1}}
\put(179,-4){\small{1}}\put(189,-4){\small{1}}
\put(100,0){\line(2,3){5}}
\put(110,0){\line(-2,3){5}}\put(120,0){\line(2,3){5}}
\put(130,0){\line(-2,3){5}}\put(140,0){\line(2,3){5}}
\put(150,0){\line(-2,3){5}}\put(160,0){\line(2,3){5}}
\put(170,0){\line(-2,3){5}}\put(180,0){\line(2,3){5}}
\put(190,0){\line(-2,3){5}} \put(124.5,3){\small{0}}
\put(144.5,3){\small{0}}\put(164.5,3){\small{0}}
\put(184.5,3){\small{0}}\put(104.5,3){\small{0}}
\put(84.5,3){\small{0}}\put(64.5,3){\small{0}}
\put(44.5,3){\small{0}}
\put(65,7.5){\line(-3,2){10}}\put(45,7.5){\line(3,2){10}}\put(85,7.5){\line(3,2){10}}
\put(105,7.5){\line(-3,2){10}}\put(125,7.5){\line(3,2){10}}
\put(145,7.5){\line(-3,2){10}}\put(165,7.5){\line(3,2){10}}
\put(185,7.5){\line(-3,2){10}}\put(94,10){\small{4}}
\put(134,10){\small{4}}\put(174,10){\small{4}}\put(54,10){\small{4}}
\put(55,14){\line(3,2){20}}
\put(95,14){\line(-3,2){20}}\put(135,14){\line(3,2){20}}
\put(175,14){\line(-3,2){20}}\put(155,23){\small{0}}
\put(75,23){\small{0}}
\put(75,27){\line(2,1){40}}\put(155,27){\line(-2,1){40}}\put(115,42){\small{5}}
\put(39,0){\line(-1,0){5}}\put(30,0){\small{4}}
\put(39,7.5){\line(-1,0){5}}\put(30,7.5){\small{3}}\put(50,14){\line(-1,0){16}}
\put(60,27){\line(-1,0){26}}\put(100,47){\line(-1,0){66}}
\put(30,14){\small{2}}\put(30,27){\small{1}}\put(30,47){\small{0}}\put(28,57)
{\small{level}}}
\end{picture}
\end{center}
\caption{A complete binary tree of height 5 corresponding to an RP
composition.}\label{fig:twintree}
\end{figure}
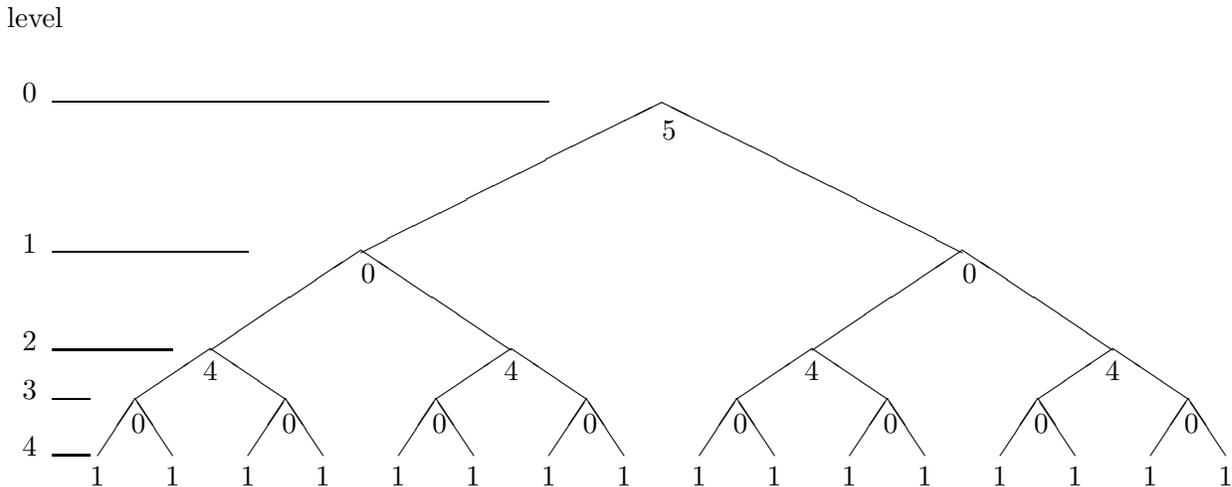
Note that the sum of labels of nodes  in such a complete binary
tree is equal to $n$. Then, reading all positive integers on such a
complete binary tree by  \textit{inorder traversal}
(left-root-right) will give an RP composition of $n$. Thus, e.g.,
the partition $\lambda=(16, 4, 4,4,4,1,1,1,1,1)$ corresponds to
the RP composition of $37$:
$37=1+1+4+1+1+1+1+4+1+1+5+1+1+4+1+1+1+1+4+1+1$.

In the same way, the bijection is easily described in general. For
an RP composition $\alpha$ of $n$, we construct a complete binary
tree $CT(\alpha)$ labeled by nonnegative integers. Denote the number
of parts of $\alpha$ by $p.$ If $\alpha=\emptyset$, then
$CT(\alpha)=\emptyset$. If $\alpha\neq \emptyset$ and $p$ is odd,
let the part in position $\frac{p-1}{2}$ be $t$. Thus
$\alpha=\alpha't\alpha'$ where $\alpha'$  is an RP composition of
$\frac{n-t}{2}.$ Now let $t$ be the root of $CT(\alpha)$ and let
$CT(\alpha')$  be the left subtree at the root. If $p$ is even, then
$\omega=\alpha'\alpha'$ where $\alpha'$
 is an RP compositions of $\frac{n}{2}$. Now let $0$ be
the root of $CT(\alpha)$, and let $CT(\alpha')$ be the left subtree
as well as the right subtree. This yields an inductive definition of
$CT(\alpha).$ It's easy to see that $CT(\alpha)$ obtained in this
way is a complete binary tree, and furthermore all nodes on the same
level have the same label. Hence it will correspond to a partition
of $n$ into powers of $2$, namely the label on level $i$ will
correspond to the multiplicity of part $2^i$. $\Box$
\section{The parity of the Catalan numbers}
It is well known that the Catalan number $C_n$ is odd iff $n=2^k-1$
for some $k$. This follows from arithmetic results of \cite{ak}, and
from bijective proofs of \cite{ed,eg,su}.

Our method of recursive palindromes provides a different, and nice
bijective proof of this fact.

On the set of binary trees of $n$ (internal) nodes we define an
involution ${\cal I}$ as follows. If $T$ is such a binary tree, we
begin at the root of $T$, and we go down the levels of $T$ until for
the first time we reach a level $L$ with the following property:
there is a node $v$ on level $L$ whose right and left subtrees are
not the same. We then interchange the right and left subtrees at
\textit{every node} on level $L$, to obtain the tree ${\cal I}(T)$.

This is evidently an involution. Its only fixed point is a complete
binary tree, and these exist only if $n=2^k-1$ for some $k$. $\Box$

We thank Professor Curtis Greene for calling this problem to our
attention.

\end{document}